\documentclass[12pt,a4paper]{amsart}

\usepackage{amsmath}
\usepackage{amscd}
\numberwithin{equation}{section}

\chardef\bslash=`\\ 

\makeatletter
\def\verbatim{\interlinepenalty\@M \@verbatim
  \leftskip\@totalleftmargin\advance\leftskip2pc
  \frenchspacing\@vobeyspaces \@xverbatim}
\makeatother
\theoremstyle{definition}

\newcounter{picture}

\newsymbol\ltimes 226E 
\newsymbol\rtimes 226F 

\newcommand{\FF}{{\mathbb F}}

\newcommand{\QQ}{{\mathbb Q}}

\newcommand{\ZZ}{{\mathbb Z}}
\newcommand{\cA}{{\mathcal A}}

\newcommand{\G}{{\Gamma}}
\newcommand{\Om}{{\Omega}}

\newcommand{\tA}{{\widetilde A}}

\newcommand{\id}{{\bf 1}} 

\begin{document}

\title[]{K-theory computations for boundary algebras of $\tA_2$ groups}

\date{August 22, 2000}
\author{Guyan Robertson }
\address{Mathematics Department, University of Newcastle, Callaghan, NSW 
2308, Australia}
\email{guyan@maths.newcastle.edu.au}
\author{Tim Steger}
\address{Istituto Di Matematica e Fisica, Universit\`a degli Studi di
Sassari, Via Vienna 2, 07100 Sassari, Italia}
\email{steger@ssmain.uniss.it}
\thanks{This research was supported by the Australian Research Council.} 
\thanks{ \hfill Typeset by  \AmS-\LaTeX}

\begin{abstract}
This is an appendix to the paper {\bf Asymptotic K-theory for groups acting on $\tA_2$ buildings}, and contains the results of the computations performed by the authors.
\end{abstract}      

\maketitle
\section{Remarks on the $\tilde A_2$ groups used in the computations} \label{computationremarks}
The results of the computations described in \cite{rsb} are given below.
We first explain some background to the terminology used in the tables.
For $q=2,3$ we give a complete list of the $\tilde A_2$ groups, following the notation of \cite{cmsz}. In these cases we indicate in the left hand column those groups that imbed as lattices in $PGL_3(K)$ for $K=\FF _2((X)),\QQ_p$ as opposed to those which are nonlinear in the sense that the associated building is not that of any linear group $PGL_3(K)$. For $q=2$ there are no nonlinear examples.
The regular and semiregular triangle presentations for $\tA_2$ groups are those described in \cite[I: Sections~4 and 5]{cmsz} respectively. Fix a projective plane over $\FF_q$.  Identify $(\FF_q)^3$ with $\FF_{q^3}$ and identify the points of the projective plane with
$\FF_{q^3}^\times / \FF_q^\times$, a cyclic group.  This cyclic group also acts on the projective plane.  The regular and semiregular triangle presentations are the triangle presentations based on this projective plane and stabilized by this cyclic group. Note that the formula given in \cite[I]{cmsz} for the number of regular/semiregular triangle presentations does not take into account equivalence of triangle presentations under
(A) other symmetries of the projective plane;
(B) passing from generators to inverse generators, and hence inverting
    the order of each triangle presentation.

 The regular and semiregular triangle presentations can also be described in the following simple way: if we number our generators $x_0,\dots,x_{q^2+q}$, a triangle presentation is regular/semiregular if it is invariant under the group $\ZZ/(q^2+q+1)$ of cyclic permutations.

If the regular/semiregular triangle presentations are defined this way, it is evident that for $(r,q^2+q+1)=1$, the permutation $j\mapsto rj$ of the
indices will take any regular/semiregular triangle presentation to another such.  To avoid duplicate groups, we have to divide through by this action.

As explained in \cite[I: Section 5]{cmsz}, a regular/semiregular triangle presentation is automatically fixed by the map $j\mapsto qj$.
Consider any triangle presentation.  If $\phi$ is an order-3  permutation of
the indices such that the triangle presentation is fixed by the map
$j\mapsto \phi(j)$, then it is easy to check that one obtains a
new triangle presentation replacing in each case $x_j x_k x_l=1$ with $y_j y_{\phi(k)} y_{\phi^2(l)}=1$. This new triangle presentation does {\it not} define the same group.  Rather, the original group $\Gamma$ can be extended to
$\Gamma^+ = \ZZ/3 \ltimes_\phi \Gamma$.
The larger group, $\Gamma^+$, acts on the building of $\Gamma$.  The
group associated to the modified triangle presentation, $\Gamma'$, corresponds to a
different index~3 subgroup of $\Gamma^+$.  Indeed, the generators of
$\Gamma'$ can be taken as $y_j=\phi\cdot x_j$.
In summary, the original group $\Gamma$ and the new group $\Gamma'$
are commensurable, each being an index~3 subgroup of a common
supergroup $\Gamma^+$.  There is a third group as well, $\Gamma''$,
whose triangle presentation is obtained by replacing in each case
$x_j x_k x_l=1$ with $z_j z_{\phi^2(k)} z_{\phi(l)}=1$.

If $q\equiv 1 \pmod 3$ then $q^2+q+1$ is a multiple of~$3$.
Consequently, one possible~$\phi$ is the translation by $(q^2+q+1)/3$.
Any regular/semiregular triangle presentation contains all triples of {\it two} of the following forms (a), (b), and (c):
\begin{itemize}
\item [(a)] $x_j x_j x_j = 1$;
\item [(b)] $x_j x_{j+(q^2+q+1)/3} x_{j+2(q^2+q+1)/3} = 1$;
\item [(c)] $x_j x_{j+2(q^2+q+1)/3} x_{j+(q^2+q+1)/3} = 1$.
\end{itemize}
If we use translation by $(q^2+q+1)/3$ as our~$\phi$ then the triangle presentations 
containing, for example, the triples (a) and (b), will be transformed
into the triangle presentations containing the triples (b) and (c).

For $q=7$ the triangle presentations obtained in this way from the Regular 
(Semiregular) presentations are called Near Regular~B and~C 
(Semiregular~B and~C). For $q=4$ there is no Semiregular presentation and
Near Regular~B and~C give isomorphic groups, called simply Near Regular. 

Another possibility for $\phi$ is the map $j\mapsto qj$.  Applied
to a regular/semiregular triangle presentation, this $\phi$ gives a triangle presentation which is 
{\it not} regular/semiregular.
For $q=2$, the regular triangle presentation, A.1 is mapped to A.2, and on a second
application of $\phi$ to A.3.  For $q=3$, the regular triangle presentation, 1.1, is
mapped to 1.2 and by a second application to 1.3.  (This is somewhat
imprecise: it is precise modulo the equivalence of triangle presentations obtained
by replacing generators with inverse generators.)  In these cases, the
abelianizations and the $K$-theory show that the groups obtained by
applying $\phi$ and by applying $\phi^2$ are not isomorphic.

Thus, for $q=4,5,7,8,9,11$, one can triple the number of groups one 
has to work with by applying the map $\phi(j)=qj$. Given a
regular/semiregular triangle presentation T, there exist new triangle presentations T$'$ and T$''$ as follows:
$$
{\text T}:\quad x_j x_k x_l = 1 \qquad\qquad
{\text T'}:\quad	    y_j y_{qk} y_{q^2l} = 1  \qquad\qquad
{\text T''}:\quad     z_j z_{q^2k} z_{ql} = 1.
$$
In particular, if the original triangle presentation was called, e.g.
$q=11$ Semiregular 6,
the two new triangle presentations are called $q=11$ Semiregular 6$'$
and $q=11$ Semiregular 6$''$, {\it even though} they are not actually semiregular.

Also included in the tables are 5-adic and 7-adic examples which were discovered recently by \textsc{H. Voskuil} \cite{v} and worked out in detail by \textsc{D. Cartwright} (private communication). In the tables we have denoted them simply ``Voskuil''.

We do not give here triangle presentations of Voskuil's groups. Likewise, we do not explain in which order we have numbered the two (seven) inequivalent semiregular triangle presentations for $q=9$ ($q=11$). A list, with labels of all the triangle presentations for which the table gives data is available as \texttt{pub/steger/triangle\_presentations.gz} by anonymous FTP from \texttt{ftp.uniss.it}. The list is also available at 

\texttt{http://maths.newcastle.edu.au/\,$\tilde{}$\,guyan/triangle\_presentations.gz}

 Regular/Semiregular triangle presentations for $q=13$ and beyond could easily be generated upon request.

\subsection{Comparison with the $K$-theory of $C_r^*(\G)$}\label{lafforgue}

Since $C_r^*(\G)$ embeds in $C(\Om)\rtimes \G$, there is a homomorphism $K_*(C_r^*(\G)) \to K_*(C(\Om)\rtimes \G)$. It is therefore worth comparing the $K$-theories of these two algebras. Let $\G_{ab}=\G/[\G,\G]$ denote the abelianization of $\G$. There is a natural homomorphism $\kappa_{\G} : \G_{ab} \to 
K_1(C_r^*(\G))$ which is rationally injective \cite{en,bv}. For comparison we have listed in the tables the abelianization of each group.
The computations suggest that barring the prime~$q=3$, the group $K_0(C(\Om)\rtimes \G)/<[\id]>$ has nonzero $q$-primary part if and only if  $\G_{ab}$ does. Away from the prime~3, the torsion part of $K_0(C(\Om)\rtimes \G)/<[\id]>$ tends to be twice $\G_{ab}$.
This fails only for the $q=2$ group B.2, for Voskuil's $q=5$ group, and for most of the $q=3$ groups.

\section{The K-group $K_0=K_1$ for $\cA(\G)=C(\Om)\rtimes \G$, $\G$ an $\tA_2$ group.}
We give below the results of the computation of $K_0$ for some small values of $q$.
In the tables $\G$ is the $\tA_2$ group, named according to \cite{cmsz} for $q=2$ and $q=3$ and as described in section \ref{computationremarks} for $q>3$.
$\G_{ab}$ is the abelianization of $\G$. 

\noindent{\bf Notation}: $[a,b,...]$ means $\ZZ_a \oplus \ZZ_b \oplus \dots$;
$m\, [a,b,...]$ means $\ZZ^m\oplus\ZZ_a \oplus \ZZ_b \oplus \dots$;
$(j)a$ means $a,a,...,a$  $j$ times. $K_0$ is the K-group of the crossed product algebra
(recall that $K_0=K_1$), and $K_0/<[\id]>$ is the K-group modulo the class $[\id]$ of the identity.

\bigskip

\centerline{
{\small
\begin{tabular}{|l|l|l|l|l|}
\hline
          & $\G$     & $\G_{ab}$    & $K_0$ & $K_0/<[\id]>$   \\ \hline
$q=2$     & A.1	  &   [(3)2,3]   &  0 [(6)2,3]      & 0 [(6)2,3]     \\		        
The $\FF _2((X))$ cases & A.1$'$  &   [(3)2,3]   &  0 [(6)2,3]   & 0 [(6)2,3]  \\
& A.2    &  [2,3,7]     &  0 [(2)2,3,(2)7] & 0 [(2)2,3,(2)7]  \\
& A.3	  &   [2,3]      &  4 [(2)2]        & 4 [(2)2]         \\
& A.4	  &    [3,9]     &  4 [3]           &   4 [3]          \\ \hline
$q=2$& B.1	  &  [3]         &  0 [3]           & 0 [3]            \\
The $\QQ _2$ cases & B.2	  &    [(2)2,3]  &  0 [(2)2,3]      & 0 [(2)2,3]       \\
& B.3	  &  [3]         &  4 []            & 4 []             \\
& C.1	  &  [(2)2,3]    &  0 [(4)2,3]      &  0 [(4)2,3]      \\ \hline
\end{tabular}
}
}
\bigskip \bigskip \bigskip
\centerline{
{\small
\begin{tabular}{|l|l|l|l|l|}
\hline
&$\G$     & $\G_{ab}$    & $K_0$ & $K_0/<[\id]>$   \\  \hline
$q=3$& 1.1	  & [(4)3]       & 26 [2] &  26  []        \\
The $\FF _3((X))$ cases& 1.1$'$	  & [(4)3]       & 26 [2] &  26  []        \\  
&   1.2  &  [(2)3,13]  & 10 [2,(2)3,(2)13]&   10 [(2)3,(2)13] \\ 
&   1.3  &  [(2)3] & 14 [2,(2)3]&   14 [(2)3] \\ 
&   1.4  &  [2,3] & 18 [(4)2]&   18 [(3)2] \\ 
&   1.5  &  [2,3,13] & 10 [(4)2,(2)13]&   10 [(3)2,(2)13] \\ 
&   1.6  &  [(2)3,9] & 14 [2,(2)3]&   14 [(2)3] \\ 
&   1.7  &  [(2)3] & 18 [2]&   18 [] \\ 
&   1.8  &  [(2)3] & 18 [2]&   18 [] \\ 
&   1.9  &  [(2)3, 9] & 14 [2,(2)3]&   14 [(2)3] \\ 
&   1.10  &  [(2)3] & 14 [2]&   14 [] \\ 
&   1.11  &  [(2)3] & 22 [2]&   22 [] \\ 
&   1.12  &  [(2)2,(2)3] & 10 [(5)2, 3]&   10 [(4)2, 3] \\ 
&   2.1  &  [2,3,13] & 10 [(4)2,13]&   10 [(3)2,13] \\ 
&   2.2  &  [2,3] & 14 [(4)2]&   14 [(3)2] \\ 
&   3.1  &  [3,13] & 10 [2,13]&   10 [13] \\ 
&   3.2  &  [3] & 14 [2]&   14 [] \\              \hline
$q=3$&   4.1  &  [4,3] & 14 [2,4]&   14 [4] \\ 
The $\QQ _3$ cases&   4.2  &  [4,3] & 14 [2,4]&   14 [4] \\ 
&   4.3  &  [2,8,3] & 10 [(2)2,8]&   10 [2,8] \\ 
&   4.4  &  [2,8,3] & 10 [(2)2,8]&   10 [2,8] \\ 
&   5.1  &  [4,3] & 14 [2,4]&   14 [4] \\ 
&   6.1  &  [(2)3] & 14 [2]&   14 [] \\ 
&   7.1  &  [(2)3] & 14 [2]&   14 [] \\ 
&   8.1  &  [2,(2)3] & 14 [(4)2]&   14 [(3)2] \\ \hline
\end{tabular}
}
}

\vfill\eject  
\centerline{
{\small
\begin{tabular}{|l|l|l|l|l|}
\hline
&$\G$     & $\G_{ab}$    & $K_0$ & $K_0/<[\id]>$   \\  \hline
$q=3$& 9.1  &  [3] & 14 [2]&   14 [] \\ 
The nonlinear&   9.2  &  [3] & 10 [2]&   10 [] \\ 
cases&   9.3  &  [3] & 14 [2]&   14 [] \\ 
&   10.1  &  [3] & 18 [2]&   18 [] \\ 
&   10.2  &  [(2)2,3] & 10 [(3)2]&   10 [(2)2] \\ 
&   10.3  & [3]  & 10 [2]&   10 [] \\ 
&   11.1  &  [3] & 14 [2]&   14 [] \\ 
&   11.2  &  [3] & 10 [2]&   10 [] \\ 
&   11.3  &  [3] & 14 [2]&   14 [] \\ 
&   12.1  &  [3,7] & 10 [2,7]&   10 [7] \\ 
&   12.2  &  [3] & 14 [2]&   14 [] \\ 
&   13.1  &  [(2)2,3] & 10 [(3)2]&   10 [(2)2] \\ 
&   13.2  &  [3] & 14 [2]&   14 [] \\ 
&   14.1  &  [3] & 14 [2]&   14 [] \\ 
&   15.1  &  [2,3] & 10 [4]&   10 [2] \\ 
&   16.1  &  [(2)3] & 14 [2]&   14 [] \\ 
&   17.1  &  [3] & 14 [2]&   14 [] \\ 
&   18.1  &  [2,3] & 14 [(2)2]&   14 [2] \\ 
&   19.1  &  [3] & 14 [2]&   14 [] \\ 
&   20.1  &  [2,3] & 10 [4]&   10 [2] \\ 
&   21.1  &  [3] & 14 [2]&   14 [] \\ 
&   22.1  &  [4,3] & 10 [2, 4]&   10 [4] \\ 
&   23.1  &  [3] & 10 [2]  &   10 [] \\ 
&   24.1  &  [3] & 10 [2]&   10 [] \\ 
&   25.1  &  [3] & 14 [2]&   14 [] \\ 
&   26.1  &  [(2)3] & 14 [2]&   14 [] \\ 
&   27.1  &  [2,3] & 14 [(2)2]&   14 [2] \\ 
&   28.1  &  [(2)3] & 14 [2]&   14 [] \\ 
&   29.1  &  [3] & 14 [2]&   14 [] \\ 
&   30.1  &  [2,3] & 10 [4]&   10 [2] \\ 
&   31.1  &  [3] & 14 [2]&   14 [] \\ 
&   32.1  &  [3] & 10 [2]&   10 [] \\ 
&   33.1  &  [3] & 14 [2]&   14 [] \\ 
&   34.1  &  [(2)2,3] & 10 [(5)2]&   10 [(4)2] \\ 
&   35.1  &  [4,3] & 10 [(3)2,4]&   10 [(2)2,4] \\ 
&   36.1  &  [3] & 14 [2]&   14 [] \\ 
&   37.1  &  [3] & 14 [2]&   14 [] \\ 
&   38.1  &  [(3)2,3] & 10 [(6)2]&   10 [(5)2]\\ 
&   39.1  &  [3] & 14 [2]&   14 [] \\ 
&   40.1  &  [8,3] & 10 [(3)2,8]&   10 [(2)2,8] \\ 
&   41.1  &  [3] & 10 [2]&   10 [] \\ 
&   42.1  &  [3] & 14 [2]&   14 [] \\ 
&   43.1  &  [3] & 14 [2]&   14 [] \\ 
&   44.1  &  [3] & 10 [2]&   10 [] \\ 
&   45.1  &  [3] & 10 [2]&   10 [] \\ 
&   46.1  & [3]  & 14 [2]&   14 [] \\ 
&   47.1  &  [3] & 14 [2]&   14 [] \\ 
&   48.1  &  [3] & 14 [2]&   14 [] \\ 
&   49.1  &  [3] & 10 [2]&   10 [] \\ 
&   50.1  &  [3] & 14 [2]&   14 [] \\ 
 \hline
\end{tabular}
}
}
\centerline{
{\small
\begin{tabular}{|l|l|l|l|l|}
\hline
&$\G$     & $\G_{ab}$    & $K_0$ & $K_0/<[\id]>$   \\  \hline
$q=3$&   51.1  &  [3] & 14 [2]&   14 [] \\ 
The nonlinear&   52.1  &  [3] & 14 [2]&   14 [] \\ 
cases&   53.1  &  [3] & 14 [2]&   14 [] \\ 
&   54.1  &  [3] & 10 [2]&   10 [] \\ 
&   55.1  &  [3] & 14 [2]&   14 [] \\ 
&   56.1  &  [3] & 10 [2]&   10 [] \\ 
&   57.1  &  [(2)2,3] & 10 [2, 4]&   10 [(2)2] \\ 
&   58.1  &  [3] & 10 [2]&   10 [] \\ 
&   59.1  & [3]  & 14 [2]&   14 [] \\ 
&   60.1  &  [3] & 14 [2]&   14 [] \\ 
&   61.1  &  [3] & 14 [2]&   14 [] \\ 
&   62.1  &  [3] & 14 [2]&   14 [] \\ 
&   63.1  &  [2,3] & 10 [(4)2]&   10 [(3)2] \\ 
&   64.1  &  [(2)2,3] & 10 [(5)2]&   10 [(4)2] \\ 
&   65.1  &  [3] & 10 [2]&   10 [] \\  \hline
\end{tabular}
}
}
\centerline{
{\small
\begin{tabular}{|l|l|l|l|l|}
\hline
&$\G$     & $\G_{ab}$    & $K_0$ & $K_0/<[\id]>$   \\  \hline
$q=4$ & Regular &  [(6)2,(2)3]  &   28 [(12)2,(6)3]   &  28 [(12)2,(6)3] \\                & Regular$'$ &[(2)2,(2)3,7]  &  28 [(4)2,(6)3,(2)7] &  28 [(4)2,(6)3,(2)7] \\    
& Regular$''$  &  [(2)2,(2)3]  & 40 [(4)2,(2)3]  &  40 [(4)2,(2)3]  \\      
& Near Regular  &	 [(2)3]  &   56 [3]  &  56 [3] \\                  
& Near Regular$'$  &  [(2)2,(2)3,7]  &  32 [(4)2,3,(2)7] &  32 [(4)2,3,(2)7] \\         
& Near Regular$''$  &  [(2)2,(2)3]  &	 32 [(4)2,3] &  32 [(4)2,3] \\ \hline             
$q=5$ & Voskuil &  [(2)2,(4)4,3]   & 62 [(14)2,(5)4,3] & 62 [(14)2,(4)4,3] \\
&Regular  &  [3,(3)5]  &  62 [4,3,(6)5]   &  62   [3,(6)5]     \\      
& Regular$'$ & [3,5,31] & 62 [4,3,(2)5,(2)31]  & 62   [3,(2)5,(2)31] \\    
&Regular$''$  &  [3,5]  &	 70 [4,(2)5]   &   70   [(2)5]  \\         
&Semiregular &  [3] & 62 [4,3]   &  62   [3]   \\         
&Semiregular$'$& [3] & 66 [4]   &    66  [] \\            
&Semiregular$''$& [3] & 66 [4]   &   66  [] \\ \hline
$q=7$ & Voskuil &   [(7)3] &   190 [2,(21)3] & 190   [(21)3] \\     
&Regular &  [(2)3,(3)7] & 190 [2,(6)3,(6)7]  & 190   [(6)3,(6)7] \\      
&Regular$'$ & [(2)3,7,19] & 190 [2,(6)3,(2)7,(2)19] &190   [(6)3,(2)7,(2)19] \\
&Regular$''$ & [(2)3,7] &	 202 [2,(2)3,(2)7]   &  202   [(2)3,(2)7]  \\ 
&Near Regular B & [(2)3] & 266 [2,3]  & 266   [3] \\                
&Near Regular B$'$& [(2)3,7,19] & 194 [2,3,(2)7,(2)19]  & 194   [3,(2)7,(2)19]  \\  
&Near Regular B$''$ & [(2)3,7]  & 206 [2,3,(2)7] & 206   [3,(2)7] \\        
&Near Regular C &	[(2)3]  &  266 [2,3]  & 266   [3]  \\              
&Near Regular C$'$ & [(2)3,7,19] & 194 [2,3,(2)7,(2)19]  & 194   [3,(2)7,(2)19] \\   
&Near Regular C$''$	& [(2)3,7] &  194 [2,3,(2)7]  & 194   [3,(2)7] \\  
&Semiregular &  [(2)3]	& 190 [2,(6)3]  &  190   [(6)3] \\       
&Semiregular$'$ &  [(2)3] &  202 [2,(2)3]  &  202   [(2)3] \\       
&Semiregular$''$  &  [(2)3] & 190 [2,(6)3]  & 190   [(6)3]  \\    
&Semiregular B  &  [(2)3] &  266 [2,3] & 266 [3] \\               
&Semiregular B$'$ &	[(2)3] & 194 [2,3]  & 194   [3] \\               
&Semiregular B$''$	& [(2)3] &  206 [2,3]  & 206   [3] \\              
&Semiregular C & [(2)3] & 266 [2,3]  & 266   [3]  \\              
&Semiregular C$'$  & [(2)3] & 194 [2,3]  & 194   [3] \\                
&Semiregular C$''$ & [(2)3] & 194 [2,3]  & 194   [3] \\ \hline 
\end{tabular}
}
}
\vfill\eject
\centerline{
{\small
\begin{tabular}{|l|l|l|l|l|}
\hline
&$\G$     & $\G_{ab}$    & $K_0$ & $K_0/<[\id]>$   \\  \hline
$q=8$ &  Regular &  [(9)2,3] & 292 [(18)2,3,7]  & 292 [(18)2,3] \\          
&Regular$'$&  [(3)2,3,73] & 292 [(6)2,3,(2)73,7] &  292 [(6)2,3,(2)73] \\    
&Regular$''$ & [(3)2,3]  & 304 [(6)2,7]   &  304 [(6)2] \\           
&Semiregular  & [3] & 292 [3,7]  &  292 [3] \\             
&Semiregular$'$ &  [3] & 296 [7]  &  296 []  \\             
&Semiregular$''$ & [3] & 300 [7]  &  300 []  \\ \hline
$q=9$ & Regular & [(7)3]  & 546 [8]  & 546  [] \\                
&Regular$'$& [(3)3,7,13]& 426 [8,(4)3,(2)7,(2)13]& 426 [(4)3,(2)7,(2)13] \\
&Regular$''$ & [(3)3] & 438 [8,(4)3]  & 438   [(4)3]  \\          
&Semiregular 1 & [3] & 546 [8] & 546  [] \\                  
&Semiregular 1$'$ & [3] & 430 [8] & 430  [] \\                  
&Semiregular 1$''$ & [3] & 434 [8] & 434  []  \\                 
&Semiregular 2 & [(4)3] & 546 [8] & 546  [] \\                 
&Semiregular 2$'$ &	[(2)3] & 434 [8,(2)3] &  [(2)3] \\            
&Semiregular 2$''$ & [(2)3,7] & 430 [8,(2)3,(2)7] & 430   [(2)3,(2)7] \\ \hline      
$q=11$ &Regular& [3,(3)11] &798 [10,3,(6)11]& 798    [3,(6)11] \\
&Regular$'$ &[3,7,11,19]& 798 [10,3,(2)7,(2)11,(2)19] & 798    [3,(2)7,(2)11,(2)19]\\
&Regular$''$&  [3,11] & 814 [10,(2)11] & 814  [(2)11] \\
&Semiregular 1& [3] & 798 [10,3] & 798    [3] \\
&Semiregular 1$'$ & [3] & 802 [10] &  802   [] \\
&Semiregular 1$''$  & [3] & 810 [10] & 810   [] \\
&Semiregular 2  &  [(3)2,3]  & 798 [10,(9)2,3] & 798  [(9)2,3] \\
&Semiregular 2$'$ &  [2,3]  &  802 [10,(3)2] &  802  [(3)2] \\
&Semiregular 2$''$  &  [2,3]  &  810 [10,(3)2] &  810 [(3)2] \\
&Semiregular 3  &  [3] &  798 [10,3] & 798  [3] \\
&Semiregular 3$'$ & [3]  & 806 [10] &  806   [] \\
&Semiregular 3$''$  & [3]  & 806 [10] & 806   [] \\
&Semiregular 4  & [3] & 798 [10,3] & 798 [3] \\
&Semiregular 4$'$ & [3] & 802 [10] &  802   [] \\
&Semiregular 4$''$  & [3] & 810 [10] & 810   [] \\
&Semiregular 5  & [3]  & 798 [10,3] & 798    [3] \\
&Semiregular 5$'$  &  [3] & 806 [10] &  806   [] \\
&Semiregular 5$''$  & [3]  &  806 [10] & 806   [] \\
&Semiregular 6   & [3]  & 798 [10,3] & 798    [3] \\
&Semiregular 6$'$  & [3] & 806 [10] & 806   [] \\
&Semiregular 6$''$  & [3] &  806 [10] & 806   [] \\
&Semiregular 7  & [3] & 798 [10,3] & 798    [3] \\
&Semiregular 7$'$  & [3] &  810 [10] &  810   [] \\
&Semiregular 7$''$  & [3] & 802 [10] & 802   [] \\
\hline
\end{tabular}
}
}


\begin{thebibliography}{RRR}



\bibitem [BV]{bv} H. Bettaieb, A. Valette, Sur le groupe $K\sb 1$ des $C^*$-alg\`ebres r\'eduites de groupes discrets. {\it  C. R. Acad. Sci. Paris S\'er. I Math.} {\bf 322} (1996), 925--928.

\bibitem[CMSZ]{cmsz} D. I. Cartwright, A. M. Mantero, T. Steger and A. Zappa, Groups acting simply transitively on the vertices of a building of type $\widetilde A_2$, I,II,\ {\it Geom. Ded.}   {\bf 47} (1993), 143--166 and 167--223.

\bibitem [EN]{en} G. Elliott and T. Natsume, A Bott periodicity map for crossed products of $C^*$-algebras by discrete groups, {\it $K$-theory} {\bf 1} (1987), 423--435. 

\bibitem[RS]{rsb} G. Robertson and T. Steger, Asymptotic K-theory for groups acting on $\tA_2$ buildings, {\it Can. J. Math.} submitted.

\bibitem[V]{v} H. Voskuil, {\it  Ultrametric Uniformization and Symmetric Spaces}, Ph.D. Thesis, Rijksuniversiteit Groningen, 1990.

\end{thebibliography}
\end{document}